\documentclass[10pt,twoside]{amsart}

\usepackage{amssymb}
\usepackage{amsxtra}
\usepackage{graphicx}
\usepackage{dcpic,pictexwd}

\theoremstyle{plain}

\newtheorem{thm}{Theorem}

\newtheorem{lemma}{Lemma}

\newcommand{\Z}{\mathbb{Z}}

\newcommand\ra{\rightarrow}
\newcommand\lra{\longrightarrow}

\newcommand{\ben}{\begin{enumerate}}
\newcommand{\een}{\end{enumerate}}
\newcommand{\disp}{\displaystyle}

\newcommand{\Oz}{P.\ Ozsv{\'a}th\,}
\newcommand{\Sz}{Z.\ Szab{\'o}\,}

\begin{document}
\title{Extending Van Cott's bounds for the $\tau$ and $s$-invariants of a satellite knot}
\author{Lawrence P. Roberts}
\maketitle

\section{Introduction}

\noindent For a knot $K \subset S^{3}$, both Khovanov homology and knot Floer homology define maps $\nu : K \lra \Z$ such that 
\ben
\item $\nu$ is a homomorphism from the group of concordance classes to $\Z$,
\item $|\nu(K)| \leq g_{4}(K)$, where $g_{4}$ is the smooth four ball genus,
\item $\nu(T_{p,q}) = \disp{\frac{(p-1)(q-1)}{2}}$ where $p, q > 0$ and $T_{p,q}$ is the $(p,q)$-torus knot.
\een 
\ \\
\noindent We assume throughout that $\nu$ is a map as above satisfying these properties. We note that these imply $\nu(\overline{K}) = -\nu(K)$, where $\overline{K}$ is the mirror of $K$, and $\nu(K \# J) = \nu(K) + \nu(J)$. Furthermore, for {\em any} orientable surface, $\Sigma$, smoothly and properly embedded in $B^{4}$, with $\partial \Sigma = K$, we have $|\nu(K)| \leq g(\Sigma)$. 
 \\
\ \\
\noindent If we denote by $K_{l,n}$, the $(l,n)$ cable of $K$, C. Van Cott
proved

\begin{thm}[Thrm. 2, \cite{VanC}]
Let $h(n) = \nu(K_{l,n}) - \disp{\frac{(l-1)n}{2}}$. Then for $n > r$, $n,r$ relatively prime to $l$, 
$$
-(l-1) \leq h(n) - h(r) \leq 0
$$
\end{thm}

\noindent In this paper we aim to generalize Van Cott's techniques to all satellites. Let $A$ be an annulus, and $P \subset A \times I$ be an embedded copy of $S^{1}$. Let $C \subset S^{3}$ be a knot. Define $S_{r}(C, P)$ to be the isotopy class of the image of $P$ under a map taking $A \times I$ to a tubular neighborhood of $C$, preserving orientations, and mapping $\partial A \times I$ to two parallel $r$-framings of $C$ (relative to the Seifert framing). Now take $A \times I \cong S^{1} \times I^{2}$, and orient $P$ so that the intersection number
of $P$ with $I^{2}$ is non-negative. Let $n_{+}$ be the number of positive intersections, and $n_{-}$ be the number of negative intersections. 

\begin{thm}\label{thrm:main}
Let 
$$
g(r) = \nu(S_{r}(C,P)) - \frac{l(l-1)}{2}r
$$ 
Then if $s > r$ and $n_{+} > n_{-}$ then
$$
-(n_{+} - 1) \leq g(s) - g(r) \leq n_{-}
$$
If $s > r$ and $n_{+} = n_{-}$ then
$$
-n_{+} \leq g(s) - g(r) \leq (n_{-} - 1)
$$
\end{thm}

\noindent Van Cott's cabling result corresponds to $n_{+} = l$, $n_{-} = 0$. $P$ can be taken to be the closure of $(\sigma_{l-1}\ldots \sigma_{1})^{m})$ for $n = l\,r + m$ where $0 \leq m < l$. We recover most\footnote{Van Cott's result applies to all values of $n$ and $r$, whereas ours applies to those congruent mod $l$} of Van Cott's result by subtracting $\frac{m(l-1)}{2}$ from $g(r)$
for each $r$ (as this is constant in $r$ it does not change the inequalities).   \\
\ \\
\section{Proof of Theorem \ref{thrm:main}}
\begin{center}
\begin{figure}
\includegraphics[scale=0.5]{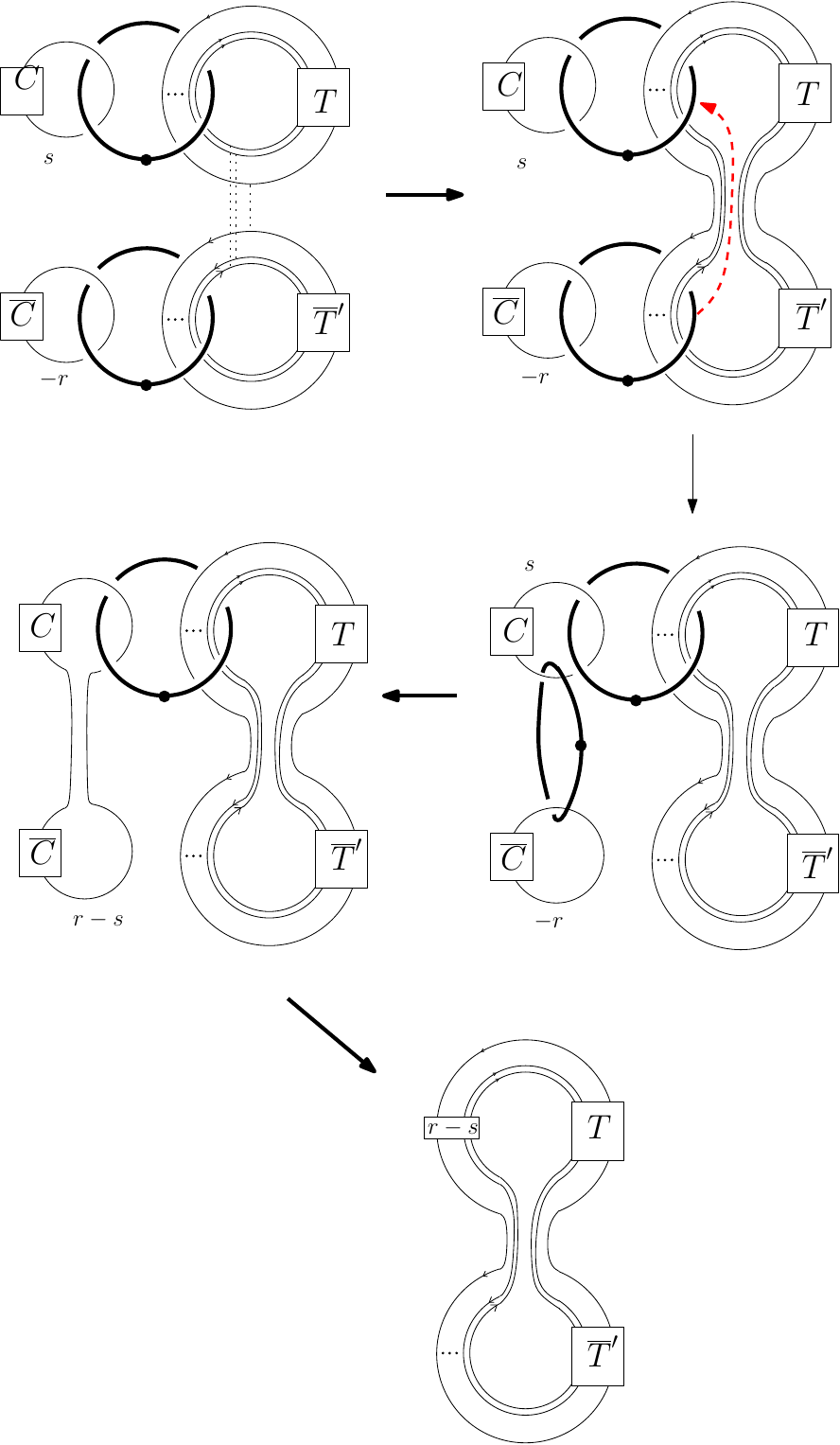}
\caption{Initial band sums and simplifications in the proof of Theorem \ref{thrm:main}. The knot $S_{-r}(\overline{C}, \overline{P})$ with tangle $\overline{T}$ was rotated by $\pi$ about the horizontal axis of the page before taking the band sums. See the next figure for an illustration. The last arrow, pointing diagonally, is a concordance of links, whereas the initial arrows are isotopies, handleslides and cancellations}\label{fig:VanCott1}
\end{figure}
\end{center}
\begin{center}
\begin{figure}
\includegraphics[scale=0.5]{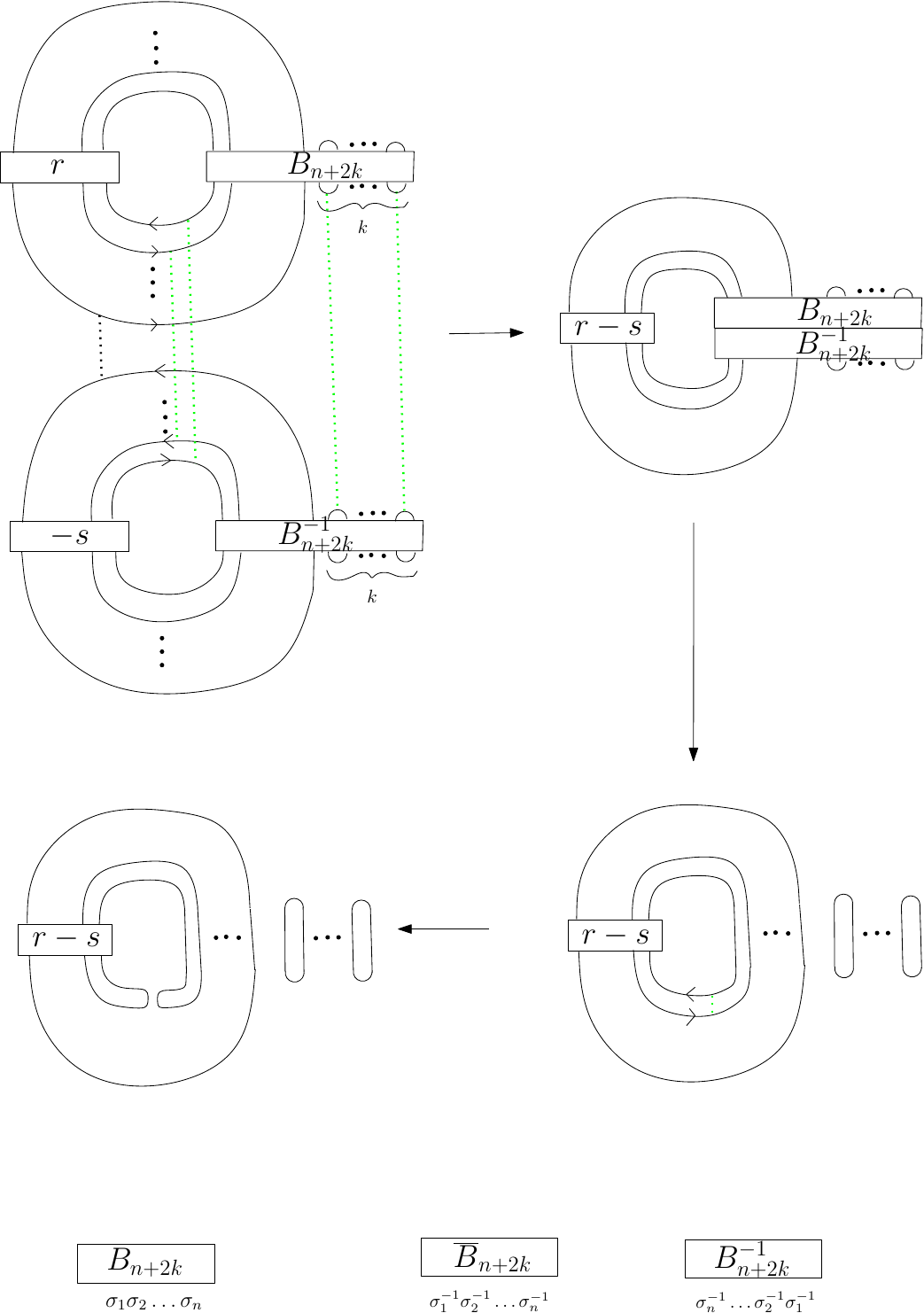}
\caption{An example beginning after we have used a concordance to remove dependence on $C$ and after we have placed $T$ in Morse position. The lowest row of diagrams shows the effect of taking the mirror image followed by a rotation of $\pi$ around the horizontal axis of the page. This is used to simplify the diagrams above. Finally, the final diagram shows the link after band summing two components of opposite orientation in a torus link with orientations altered. The ``cup'' can be isotoped around the diagram to form an unknot with the ``cap''.}\label{fig:VanCott2}
\end{figure}
\end{center}
\noindent We start with $P \subset S^{1} \times I^{2}$, as above. We may isotope $P$ so that $\{p\} \times I^{2}$ intersects $P$ transversely in $n_{+}$ positive intersection points, and $n_{-}$ negative intersection points. There is an ambient isotopy of $S^{1} \times I^{2}$, which preserves the framing, but ensures that $P$ is a product on $(p-\pi/2, p+ \pi/2)$. Consequently we may assume that $P$ has the following schematic representation:

\begin{center}
\includegraphics[scale=0.5]{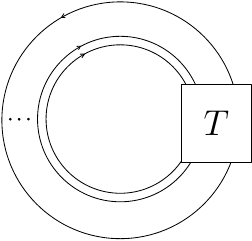}
\end{center}

\noindent Within the box, $P$ induces a tangle, $T_{P}$, which can be isotoped so that the projection $pr_{1} : (p+\pi/2, p-\pi/2) \times I^{2} \ra (p+\pi/2, p-\pi/2)$ has only local maxima and minima for critical points. We can then isotope this portion of $P$ to have the following structure for some braid on $n_{+} + n_{-} + 2k$ strands:

\begin{center}
\includegraphics[scale=0.5]{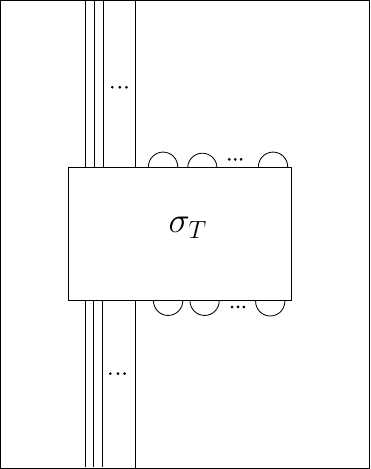}
\end{center}

\noindent This is a Morse presentation of the tangle, and can be obtained by ``pulling'' the maxima up and to the right, and the minima down and to the right, in any manner that remains embedded. Choose such a representation for $P$, and call the braid $\sigma_{T}$. Note that the isotopies induce on the boundary a map that is isotopic to the identity (indeed, this could have been done rel boundary if we had preferred), so we may use this presentation of $P$ in the construction of a satellite without affecting the isotopy class of the satellite.\\
\ \\
\noindent We will think of $S_{r}(C,P)$ as the image of the knot $P$ under the diffeomorphism of $S^{3}$ induced by canceling the $1$-handle/$2$-handle pair depicted in the following surgery diagram for $B^{4}$: 

\begin{center}
\includegraphics[scale=0.5]{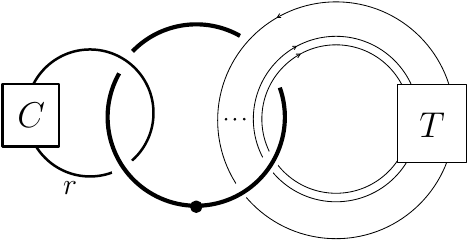}
\end{center}

\noindent Note that the $r$ on $C$ provides the twisting defined above, and the meridian of the $1$-handle corresponds to the framing on $S^{1} \times D^{2}$. To see that we obtain $S_{r}(C,P)$ under the diffeomorphism, slide $P$ over $C$ $n_{+} + n_{-}$ times so that it no longer runs over the $1$-handle. \\
\ \\
\noindent Consider $K_{s} = S_{s}(C,P)$ and $K_{r} = S_{r}(C,P)$ with $s > r$. Let $K'_{-r}$ be the mirror of $K_{r}$ with reversed strand orientation. We will start by considering $K_{s} \# K'_{r}$, see as an example Figure \ref{fig:VanCott1}. Let $n= n_{+} + n_{-}$ be the number of strands, where $n_{+}$ are oriented up and $n_{-}$ are oriented down, and band sum $n - 1$ times as in the figure. The result is a satellite link of 
$C \# \overline{C}$ with twisting number $s-r$. As such it is concordant to a satellite with companion the unknot and the same pattern, with twisting number $s-r$ after employing

\begin{lemma}\label{lem:satcon}
Let $\overline{C}$ be the mirror of $C$. The knot $S_{m}(\overline{C}, S_{r}(C,P))$ is concordant to $S_{m + r}(U, P)$ where $U$ is the unknot.
\end{lemma}

\noindent {\bf Proof:} Since $C \# \overline{C}$ is a ribbon knot, we can find a slice disk for it. If we trivialize a neighborhood of this disc to obtain a region of $B^{4}$ diffeomorphic to $D^{2} \times D^{2}$. If we take $D^{2} \times \{0\}$ to be the slice disc, we can construct $n$ parallel copies by choosing $n$ points, $x_{i}$, in the second factor, and taking the image of $D^{2} \times \{x_{i}\}$ under the diffeomorphism for each $i$. In $S^{3}$ this gives a link formed by $m+r$ parallel copies of $C \# \overline{C}$. Each copy is a longitude since it bounds a disc disjoint from the slice disc. We place this configuration close to $S_{m+r}(U,P)$, and orient the longitudes in such a way that we can perform $n$ band sums and obtain an oriented knot. This knot is the same as $S_{m}(\overline{C}, S_{r}(C,P))$. $\Diamond$ \\
\ \\
\noindent We note that the argument in the lemma does not require $P$ to be a knot, provided concordance for links is the relation generated by collections of disjoint, smoothly embedded annuli in $S^{3} \times I$. \\
\ \\
\noindent The resulting link, $L$, is composed of two pieces, a region of $s - r$ full twists, and the composition of two tangles on $n$ strands. If $T$ is one tangle, then the other is $\overline{T}$, the mirror with the orientations on all components reversed, {\em and} the ends of $D^{2} \times I$ switched, i.e. the tangle is turned upside down. As above we consider $T$ in Morse position, i.e. a composition of a trivial $n$ tangle with some number of cups, a braid $\sigma_{T}$, the same number of caps, followed by a trivial $n$-tangle. The mirror tangle is then given by the same diagram with each generator in $\sigma$ replaced by its inverse. If we change orientations, rotate $180^{\circ}$ about the horizontal axis, and then compose, we are composing $T$ with a tangle whose Morse position is determined by $(\sigma_{T})^{-1}$. We can
then band sum the caps of $T$ with the cups of $\overline{T}$ to obtain a diagram isotopic to a trivial $n$ tangle and some number of circles (since $\sigma_{T} (\sigma_{T})^{-1} = I$). The circles can then be made to bound discs disjoint from one another. We have the same number of discs
as bands, thus the the corresponding cobordism of tangles can be applied locally to $L$, see Figure \ref{fig:VanCott2}.\\
\ \\
\noindent The resulting link, $L'$, has the same projection as $T_{n,n(s-r)}$ but with different orientations on the strands. Now find two adjacent strands with opposite orientations and band connect them. The result is a new link with one unlinked unknot. This can be seen by noting that $L'$
is the same as taking $s-r$\  $-1$-framed unknots parallel to the axis of the closure of the trivial $n$ braid, with appropriate orientations. The braid sits in a copy of $S^{1}\times D^{2}$ in which there is an oriented annulus bounded by the two strands. The image of this annulus under blowing down the $-1$-components provides the required degeneration. If $n_{+} \geq n_{-}$ we can repeat this to pair off all the $n_{-}$ strands\\
\ \\
\noindent {\bf Case i:} If $n_{+} = n_{-}$, repeating this pairing will ultimately consume all the strands. We will thus
have constructed a surface in $B^{4}$ with boundary the knot $K_{s} \# K'_{-r}$, and used $n-1 + n_{+} + \#\mathrm{caps}(T)$ bands and $\#\mathrm{caps}(T) + n_{+}$ discs in the process. Consequently, the surface has Euler characteristic $1 - n = 1 - 2n_{+}$. Therefore, the genus of this surface, which has one boundary, is $g=n_{+} = n_{-}$.  Since $|\nu(K)| \leq g_{4}(K)$,  and $\nu(K'_{r}) = -\nu(K_{r})$ we have
$$
\big|\nu(K_{s}) - \nu(K_{r}) \big| \leq n_{+}
$$
Since $l=0$, $g(r) = \nu(K_{r})$, so we can conclude that $- n_{+} \leq g(s) - g(r) \leq n_{-}$.  \\
\ \\
\noindent However, we may construct a different surface. There are at least two strands intersecting $I^{2}$ in this case. Pair all but two, oppositely oriented strands, as before.  This introduces $n_{+} - 1$ bands and $n_{+} - 1$ discs. The remaining two strands bound an annulus, twisted positively $s-r$ times. We may now band the two strands together as in:  

\begin{center}
\includegraphics[scale=0.3]{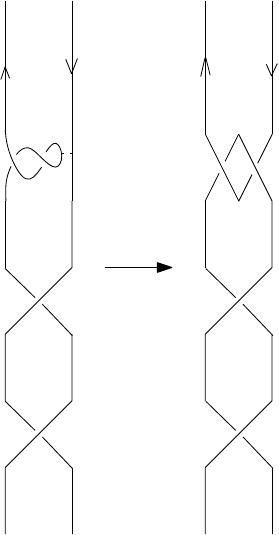}
\end{center}

\noindent The resulting knot is the $s-r$ twisted negatively clasped Whitehead double of the unknot, $D_{-}(U, s-r)$ in the notation of \cite{Hed2}. However, $\nu(D_{-}(U, s-r)) = -1$ when $s-r > 0$, since by Theorem 2 of \cite{LivN}

\begin{lemma}
For each $K$, $\nu(D_{+}(K,t)) = 1$ for $t \leq \mathrm{TB}(K)$. 
\end{lemma}

\noindent In particular $\overline{D_{-}(U, s-r)} = D_{+}(U, r-s)$, and $r-s \leq -1 = \mathrm{TB}(U)$. Consequently, $\nu(\overline{D_{-}(U, s-r)}) = 1$ and $\nu(D_{-}(U, s-r)) = -1$. We have added another band in order to do this. Consequently, we have used $n-1 + \#\mathrm{caps}(T) + n_{+} - 1 + 1$ bands and $\#\mathrm{caps}(T) + n_{+} - 1$ to create a surface with two boundaries: $K_{s} \# K_{r}'$ and $D_{-}(U, s- r)$. This surface has
Euler characteristic $n = 2n_{+}$ and genus $n_{+}$ as before. Therefore, the knot $K_{s} \# K_{r}' \# D_{+}(U, r-s)$ bounds a surface of genus
$n_{+}$ and 
$$
\big|\nu(K_{s}) - \nu(K_{r})  + 1\big| \leq n_{+}
$$
Consequently, $-(n_{+} + 1) \leq g(s) - g(r) \leq n_{-} - 1$. As both sets of inequalities must be true, we have
$ - n_{+} \leq g(s) - g(r) \leq n_{-} - 1$.  \\
\ \\
\noindent {\bf Case ii:} When $l = 1$, i.e. $n_{+} = n_{-} + 1$, we can pair the strands as above to remove all the strands contributing to $n_{-}$. This requires $n_{-}$ band sums, and $n_{-}$ disc attachments. This results in an unknot, which we can fill with a disc. The surface so constructed has $n-1 + \#(caps) + n_{-}$ bands and $\#(caps) + n_{-} + 1$ discs, for an Euler characteristic
of $2 - n = 2 - (n_{+} + n_{-}) = 2 - (n_{+} + n_{+} - 1) = 3 - 2n_{+}$. The surface's genus, therefore, is $n_{+} - 1$. Consequently,
$$
\big|\nu(K_{s}) - \nu(K_{r}) \big| \leq (n_{+} - 1)
$$
When $l=1$, $g(r) = \nu(K_{R})$ so this implies $- (n_{+} - 1) \leq g(s) - g(r) \leq n_{+} - 1$. Since $l = 1$, $n_{+} - 1 = n_{-}$, so we conclude that $- (n_{+} - 1) \leq g(s) - g(r) \leq n_{-}$.\\
\ \\
\noindent {\bf Case iii:} When $l > 1$, we can again band the strands contributing to $n_{-}$ to some of the strands contributing to $n_{+}$. This gives $n_{-}$ new bands, and $n_{-}$ unknots which we fill with discs. The remaining link is the torus link $T_{l,l(s-r)}$. 
We may now use $n_{+} - n_{-} - 1$ additional band sums to convert to either one of the torus knots $T_{l,l(s-r) + 1}$ or $T_{l,l(s-r) - 1}$ (see Figure 2 of \cite{VanC}). Since $s > r$ and $l > 1$, $l(s-r) \pm 1 > 0$.  \\

\begin{center}
\includegraphics[scale=0.3]{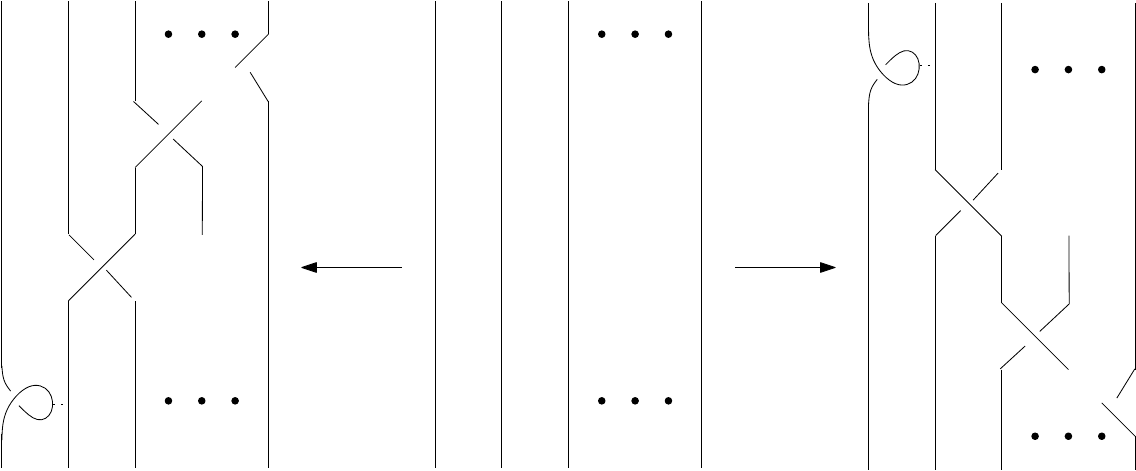}
\end{center}

\noindent We have used a total of $n - 2 +  n_{+} - n_{-}  +  n_{-} + \#(caps)$ band sums when $l > 1$. On the other hand, in this case we also fill $\#\mathrm{caps}(T) + n_{-}$ unknots by discs. The surface we have constructed has Euler characteristic $-n+2 - n_{+} + n_{-} = -2n_{+} + 2$ and two boundary components. Therefore, its genus is $g = n_{+} - 1$. We can conclude that both $K_{s} \# K'_{-r} \# \overline{T_{l, l(s-r) \pm 1}}$ bound smoothly embedded orientable surfaces with genus $g= n_{+} - 1$.\\
\ \\
\noindent Consequently
$$
\big| \nu(K_{s} \# K'_{-r} \# \overline{T_{l, l(s-r) \pm 1})} \big| \leq g
$$
Using the properties of $\nu$ listed in the introduction, we may write down two inequalities:
$$
\big| \tau(K_{s}) - \tau(K_{r}) - \frac{(l-1)l(s-r)}{2} \big| \leq  g
$$
and
$$
\big| \tau(K_{s}) - \tau(K_{r}) -  \frac{(l-1)l(s-r)}{2} + (l-1) | \leq g 
$$
Let $g(r) = \tau(K_{r}) - \frac{l(l-1)}{2}r$, then the first inequality is equivalent to $-g \leq g(s) - g(r) \leq g$ whenever $s \geq r$. This simplifies to $-(n_{+} - 1) \leq g(s) - g(r) \leq (n_{+} - 1)$. The second inequality is $-g-l+1\leq g(s) - g(r) \leq g - l +1$. As $g = n_{+} - 1$, $g - l + 1 = n_{-}$ and $-g-l+1 = -(n_{+} - 1)  - (n_{+} - n_{-}) + 1 = $ $-2n_{+} + n_{-} + 2$, so this simplifies to $2 - 2n_{+} + n_{-} \leq g(s) - g(r) \leq n_{-}$.  Both sets of inequalities must be true, but when $l > 1$, $n_{+} - 1 > n_{-}$, and $- (n_{+} - 1) \geq  -2(n_{+}-1) + n_{-}$, the stricter bounds on each side yield $-(n_{+} - 1) \leq g(s) - g(r) \leq n_{-}$.  $\Diamond$\\
\ \\

\end{document}